\documentclass[preprint,1p]{elsarticle}
\usepackage{graphicx}
\usepackage{amssymb}
\usepackage{amsmath,amsthm}
\usepackage{amsfonts}
\usepackage{graphicx}
\usepackage{graphics}
\usepackage{color}
\usepackage{latexsym}
\usepackage{url}
\def\biglf{\par\bigskip\noindent}
\def\fs#1{\footnotesize#1}
\newcommand{\R}{\ensuremath{\mathbb{R}}}

\newcommand{\calp}{\ensuremath{\mathcal{P}}}

\def\bfm#1{\protect{\makebox{\boldmath $#1$}}}

\def\b {\bfm{b}}

\def\u {\bfm{u}}

\def\p {\bfm{p}}

\def\fa {\mbox{ for all }}

\def\cala{\mathcal{A}}
\def\calb{\mathcal{B}}

\newcommand{\bql}[1]{%
\begin{equation}\label{#1}%
}
\def\RSlabel#1{\label{#1}%
}

\newtheorem{Theorem}{Theorem}[section]

\theoremstyle{definition}
\newtheorem{Definition}[Theorem]{Definition}
\theoremstyle{remark}

\numberwithin{equation}{section}

\begin{document}
\begin{frontmatter}
\title{Direct meshless local Petrov-Galerkin (DMLPG) method:\\
A generalized MLS approximation}%
\author[DM]{D. Mirzaei}
\ead{d.mirzaei@sci.ui.ac.ir}%
\address[DM]{Department of Mathematics, University of Isfahan, 81745-163 Isfahan, Iran. }
\author[RS]{Robert Schaback\corref{corRS}}
\ead{schaback@math.uni-goettingen.de}%
\cortext[corRS]{Corresponding author}
\address[RS]{Institut f\"ur Numerische und Angewandte Mathematik,
Universit\"at G\"ottingen, Lotzestra\ss{}e 16-18,
D--37073 G\"ottingen, Germany}%


\begin{keyword}
Meshless Methods \sep
Moving least squares \sep
Meshless Local Petrov Galerkin methods \sep
MLPG  \sep
shape functions \sep
diffuse derivatives
\end{keyword}


\begin{abstract}
The {\em Meshless Local Petrov--Galerkin} (MLPG) method
 is one of the popular meshless methods
that has been used very successfully to solve several types of boundary
value problems since the late nineties. In this paper, using a generalized
moving least squares (GMLS) approximation, a new {\em direct}
MLPG technique, called DMLPG,
is presented. Following the principle
of meshless methods to express
everything ``entirely in terms of nodes'',
the generalized MLS recovers test functionals directly from values at nodes,
without any detour via shape functions.
This leads to a cheaper and even more accurate
scheme. In particular, the complete absence of shape functions
allows numerical integrations in the weak forms of the problem to be done
over low--degree polynomials instead of complicated shape functions.
Hence, the standard MLS shape
function subroutines are not called at all.
Numerical examples illustrate the superiority
of the  new technique over the classical
MLPG. On the theoretical side,
this paper
discusses
stability and convergence for the new discretizations
that replace those of the standard MLPG. However, it does not
treat stability, convergence, or error estimation for the MLPG as a whole.
This should be taken from the literature on MLPG.
%
\end{abstract}
\end{frontmatter}
\section{Introduction}

The {\em Moving Least Squares method} (MLS)
was introduced as an approximation technique by
Lancaster and Salkauskas \cite{lancaster-salkauskas:1981-1},
inspired by the pioneering work of Shepard
\cite{shepard:1968-1} and McLain \cite{mclain:1974-1, mclain:1976-1}.
Since the numerical approximations of MLS are based on a cluster of
scattered nodes instead of interpolation on elements,
many meshless methods for the numerical
solution of differential equations were based on the MLS method in recent years.
As an important example of
such methods, we mention the {\em Meshless Local Petrov-Galerkin} (MLPG)
method introduced by S.N.
Atluri and his colleagues
\cite{atluri:2005-1,atluri-shen:2002-book,atluri-zhu:1998-1}.
It is a {\em truly} meshless method in weak form
which is based on local subdomains, rather than a single global domain.
It requires neither domain elements nor background cells in either the
approximation
or the integration.

In MLPG and other MLS based methods, the stiffness matrix
is provided by integrating over MLS shape
functions or their derivatives. These shape functions are
complicated and have no closed forms.
To get accurate results, numerical quadrature with many integration points is required.
Thus the MLS subroutines must be
called very often, leading to high computational costs.
In contrast to this,
the stiffness matrix in finite element methods (FEMs) is
constructed by integrating over {\em  polynomial} basis functions
which are much cheaper to evaluate. This relaxes the cost of
numerical integrations somewhat. For an account of the importance of
numerical integration within meshless methods, we refer the reader to \cite{babuska-et-al:2009-1}.

This paper avoids integration over MLS shape functions in MLPG and replaces it by
the much cheaper integration over polynomials. It ignores shape functions completely.
Altogether, the method is simpler, faster and more accurate than the original MLPG method.
We use a generalized form of the MLS which directly
approximates boundary conditions
and local weak forms, shifting the numerical integration
into the MLS itself, rather than into an outside loop
over calls to MLS routines. We call this approach
{\em Direct Meshless Local Petrov-Galerkin} (DMLPG) method.
The convergence rate of MLPG and DMLPG seems to be the same,
but thanks to the simplified computation,
the results of DMLPG often are more precise than
the results of MLPG. All of this is confirmed by numerical examples.
\section{Meshless methods}\RSlabel{SecMM}
Whatever the given problem is,
{\em meshless methods}
construct solutions from a {\em trial space} $U$ whose functions are
parametrized ``{\it entirely in terms of nodes}''
\cite{belytschko-et-al:1996-1}. We let these nodes
form a set $X:=\{x_1,\ldots,x_N\}$. Then the functions $u$
of the linear trial space $U$ are parametrizable by
their values on $X$ iff the linear functionals
$\delta_{x_1},\ldots,\delta_{x_N}$ are linearly
independent on $U$. This implies that there must be
some basis $u_1,\ldots,u_N$ of $U$ such that the $N\times N$ matrix of values
$u_j(x_k)$ is invertible, but we are not interested in knowing or
constructing this basis.
We only assume that the discretized problem is set up with
a vector
$$
\u=(u(x_1),\ldots,u(x_N))^T
$$
of unknowns
in ``meshless'' style, and all data have to be expressed in terms of these.

Furthermore, we assume the discretized problem to consist of equations
\bql{eqlkubk}
\lambda_k(u)=\beta_k,\;1\leq k\leq M,
\end{equation}
where we have $M\geq N$ linear functionals $\lambda_1,\ldots,\lambda_M$ and $M$
prescribed real values $\beta_1,\ldots,\beta_M$. Section
\ref{Sectweakform} will describe how this is done for standard linear PDE
problems, including the variations of the MLPG.

The upshot of all meshless methods now is to provide good estimates
$\widehat \lambda_k$ of all $\lambda_k$ using only values at nodes. This means
that one has to find real numbers $a_j(\lambda_k)$ with
\bql{eqlkajkuj}
\widehat{\lambda_k(u)}=\displaystyle{\sum_{j=1}^Ma_j(\lambda_k)u(x_j) \approx
  \lambda_k(u)   }
\fa k,\;1\leq k\leq M.
\end{equation}
Putting the $a_j(\lambda_k)$ into an $M\times N$ matrix
$\cala$, one has to solve the possibly overdetermined linear system
\bql{eqAub}
\cala\,\u=\b
\end{equation}
with
$\b=(\beta_1,\ldots,\beta_M)^T$.

Note that we do not mention {\em shape functions} at all.
They are not needed to set up a linear system
in terms of values at nodes. The goal just is to
find good estimates for the target functionals $\lambda_k$ in terms of the
values at nodes, e.g. via (\ref{eqlkajkuj}), to set up the matrix $\cala$.
Note that in some cases, e.g. when the functionals $\lambda_k$
are derivatives at points, this is just a variation of a finite--difference approach.

In a second stage, users might want to
evaluate $u$ at other places than in the nodes $x_j$. This is
a problem of recovery of functions from discrete data values, and
completely independent of PDE solving. There are various possibilities
to do so, including the standard MLS with its shape functions,
but we do not comment on these techniques here.
\section{Generalized moving least squares (GMLS) approximation}\RSlabel{SecMLS}
Before we show how to discretize PDEs in the form
(\ref{eqlkubk}), we focus on how to find
good estimates of functional values $\lambda(u)$ in terms of nodal values
$u(x_1),\ldots,u(x_N)$.
The classical MLS approximates $\lambda(u)=u(x)$ from nodal values,
minimizing a certain weighted discrete $l_2$ norm.
But in view of the previous section, we need
more general functionals. Therefore we employ a generalized version of Moving Least Squares,
adapted from \cite{mirzaei-et-al:2012-1}.

Let $u\in C^m(\Omega)$ for some $m\geq 0$,
and let $\{\mu_j(u)\}_{j=1}^N$
be a set of continuous linear functionals
$\mu_j$ from the dual $C^m(\Omega)^*$ of $C^m(\Omega)$.
For a fixed given functional $\lambda\in C^m(\Omega)^*$,
our problem is the approximate recovery of the
value $\lambda(u)$ from the values $\{\mu_j(u)\}_{j=1}^N$.
This fits into the preceding section for $\lambda=\lambda_k,\,1\leq k\leq M$
and $\mu_j(u)=u(x_j),\,1\leq j\leq N$.

The functionals $\lambda$
and $\mu_j$, $1\leq j\leq N$, can, for instance, describe point evaluations
of $u$, its
derivatives up to order $m$, or some local integrals that contain $u$ or its
derivatives in their integrands. In particular, we shall use functionals of the
form (\ref{localweakformeq}) arising from local weak forms, or simple
point evaluation functionals on the Dirichlet part of the boundary.
\biglf
The approximation $\widehat{\lambda(u)}$ of $\lambda(u)$ should be a linear
function
of the data $\mu_j(u)$, i.e. it should have the form
\bql{eq1}
\widehat{\lambda(u)}
=\sum_{j=1}^Na_j(\lambda)\mu_j(u),
\end{equation}
and the coefficients $a_j$ should be linear in $\lambda$.
We already saw this in (\ref{eqlkajkuj}).
As in the classical MLS, we assume the approximation equation
(\ref{eq1}) to be exact for a finite dimensional subspace
$\calp=\mathrm{span}\{p_1,p_2, \ldots , p_Q\}\subset C^m(\Omega)$, i.e.
\bql{eq2}
\sum_{j=1}^Na_j(\lambda)\mu_j(p)=\lambda(p) \fa p\in \calp.
\end{equation}
As in the classical MLS, we employ the standard technique of minimizing
\bql{eq5}
\frac{1}{2}\sum_{j=1}^Na_j^2(\lambda)/w_j
\end{equation}
as a function of the coefficients $a_j(\lambda)$ subject to the linear
constraints (\ref{eq2}),
where we use positive weights $w_1,\ldots,w_N$ which later
will be chosen in a specific way to localize the approximation,
provided that $\lambda$ is a functional acting locally.
Anyway, the weights should depend on the  functionals $\lambda$ and $\mu_j$.
In most cases, the functional $\lambda$ will be localized at
some point $x$, and then we shall use the standard MLS weights for evaluation at $x$.

The GMLS approximation $\widehat{\lambda(u)}$
to $\lambda(u)$ can also be obtained
 as $\widehat{\lambda(u)}=\lambda(p^*)$,
where $p^*\in \calp$ is minimizing
the weighted least-squares error functional
\bql{eq3}
\sum_{j=1}^N\big(\mu_j(u)-\mu_j(p)\big)^2 w_j,
\end{equation}
among all $p\in \calp$.
This problem is independent of the functional $\lambda$
and can be efficiently applied for several functionals $\lambda$
for fixed functionals $\mu_j$. This may simplify certain calculations a lot,
provided that several functionals have to be estimated based on the same local data.
Details are in \cite{mirzaei-et-al:2012-1}, including error bounds
for the recovery and a proof that (\ref{eq1}) holds for
$\widehat{\lambda(u)}=\lambda(p^*)$ with the optimal solution $p^*$
of (\ref{eq3}) and the optimal solution $a_j^*(\lambda)$
of (\ref{eq5}).
However, in meshless methods, we need more than
the single value
$\widehat{\lambda(u)}=\lambda(p^*)$,
since we finally need the solution as a
vector $a^*(\lambda)\in\R^N$ with $N$ values $a_j^*(\lambda),\,1\leq j\leq N$.

In the special case
\bql{funcdef1}
\lambda(u)=(\delta_{x} D^\alpha)(u),
\end{equation}
the derivatives of $u$ of order $|\alpha|$ are recovered.  They are called
 {\em GMLS approximation derivatives} in \cite{mirzaei-et-al:2012-1}.
Some authors call them {\em diffuse} derivatives, but they not ``diffuse'' in
any way. They are very good direct recoveries of the derivatives of $u$, but not
coincident with the corresponding derivatives of the shape functions of
the classical MLS solution. Our GMLS approach does not even have shape
functions. Instead, derivatives are estimated directly from nodal values,
avoiding the inefficient detour via classical derivatives of shape functions.

Note that the use of polynomials is not mandatory, and the resulting
values $a_j(\lambda)$ will be independent of the chosen basis of $\mathcal{P}$.
However, choosing a good basis of $\mathcal{P}$ will improve stability, and
the following discussion shows that $\mathcal{P}$ should
have the property that $\lambda(p)$ should be easy to evaluate for $p\in \mathcal{P}$.

Even if a different numerical method is used to minimize
(\ref{eq5}) or (\ref{eq3}),
the
optimal solution $a^*(\lambda)\in\R^N$ can be written as
\bql{aWl}
{a^*(\lambda)=WP^T(P\,W\,P^T)^{-1}\lambda(\p)}  
\end{equation}
where $W$ is the diagonal matrix carrying the weights $w_j$ on its diagonal,
{$P$} 
is the $N\times Q$ matrix of values $\mu_j(p_k),\;1\leq j\leq N,\,1\leq
k\leq Q$,
and 
{$\lambda(\p)\in\R^Q$} 
is the vector with values $\lambda(p_1),\ldots,\lambda(p_Q)$
of $\lambda$ on the basis of $\calp$.
Thus it suffices to evaluate $\lambda$ on the space $\calp$, not on
a certain trial space spanned  by certain shape functions.
This will significantly speed up numerical calculations, if the functional
$\lambda$ is complicated, e.g. a numerical integration against a test function.
Standard examples are functionals of the form
$$
\lambda(u)=\int_K w(x) \,Lu(x)dx
$$
where $L$ is a linear differential operator preserving polynomials
or just the identity,
and $w$ is some polynomial test or weight function. Such functionals will arise
for PDE problems in weak form in the next section. Then
our generalized MLS technique will perform integration only over
polynomials, if we use polynomials as the space $\mathcal{P}$.
Note that this generalizes to any type of functional:
we finally just have to evaluate it on a polynomial. No other calls to
MLS routines are necessary, because we do not apply the functional to shape
functions.
\section{Problems in local weak forms}\RSlabel{Sectweakform}
We now write linear PDE problems in the discretized form (\ref{eqlkubk}),
with special emphasis on the Meshless Local Petrov Galerkin Method.

Although the technique proposed in this paper can be used for a wide class of PDEs,
we illustrate our approach for the Poisson problem
\bql{poissoneq}
\begin{array}{rcll}
     \Delta u(x) &=& f(x),           & x\in \Omega, \\
     u(x)        &=& u_D(x),         & x\in\Gamma_D,\\
     \frac{\partial u}{\partial n}(x)&=& u_N(x), & x\in\Gamma_N
\end{array}
\end{equation}
where $f$ is a given source function,
the bounded
domain $\Omega\subset\R^d$ is enclosed by
the boundary
$\Gamma=\Gamma_D\cup\Gamma_N$, $u_D$ and $u_N$ are the prescribed
Dirichlet and Neumann data, respectively, on the Dirichlet boundary $\Gamma_D$ and on the
Neumann boundary $\Gamma_N$, while $n$ is the outward normal direction.

The simplest way of discretizing the problem
in the form (\ref{eqlkubk}) is direct and global collocation.
In addition to the trial nodes $x_1,\ldots,x_N$ for
obtaining nodal solution values,  we can choose finite point sets
$$
Y_\Omega\subset\Omega,\;Y_D\subset
\Gamma_D,\;Y_N\subset\Gamma_N,\;Y:=Y_\Omega\cup Y_D\cup Y_N,\;|Y|=M
$$
and discretize the problem by $M$ functionals
\bql{eqfullcoll}
\begin{array}{rcrcll}
     \lambda_i(u)&=&\Delta u(z_i) &=& f(z_i),           & z_i\in Y_\Omega\subset \Omega, \\
     \lambda_j(u)&=&u(z_j)        &=& u_D(z_j),         & z_j\in Y_D\subset \Gamma_D,\\
     \lambda_k(u)&=&\frac{\partial u}{\partial n}(z_k)&=& u_N(z_k), & z_k\in Y_N\subset \Gamma_N
\end{array}
\end{equation}
using some proper indexing scheme. In MLPG categories, this is MLPG2 \cite{atluri:2005-1,atluri-shen:2002-book}. All
functionals are local, and {\em strong} in the sense that they do not involve
integration over test functions.

For FEM--style global weak discretization, one can keep the second and third
part of (\ref{eqfullcoll}), but the first can be weakened using the Divergence
Theorem.
With sufficiently smooth test functions $v_i$, we get
$$
\lambda_i(u):=\int_{\Gamma}(\nabla u\cdot n) v_i \ d\Gamma - \int_{\Omega}
\nabla u\cdot \nabla v_i \ d\Omega =
\int_{\Omega} fv_i\ d\Omega
$$
as a replacement of the first functionals in (\ref{eqfullcoll}),
leading again to (\ref{eqlkubk}).

Following the original MLPG method, instead of transforming (\ref{poissoneq})
into a global weak form,  we construct weak forms over local subdomains
$\Omega_\sigma^y$ which are small regions taken around nodes $y$
in the global domain $\Omega$. The local subdomains could theoretically
be of any geometric shape and size. But for simplicity they
are often taken to be balls $B(y,\sigma)$ intersected with $\Omega$ and
centered at $y$ with radius $\sigma$,
or squares in 2D or cubes in 3D centered
at $y$ with sidelength $\sigma$,
denoted by $S(y,\sigma)\cap\Omega$. The variable
$\sigma$ parametrizes the local subdomain's size, and we denote the boundary
within $\Omega$ by $\Gamma_\sigma^y:=\Omega\cap \partial \Omega_\sigma^y$.
We call a node $y$ {\em internal}, if the
boundary $\partial \Omega_\sigma^y$ of the local subdomain $\Omega_\sigma^y$ does
not intersect $\Gamma$.

The derivation of the local weak form starts with the local integral
\bql{integralformeq}
\int_{\Omega_\sigma^y}\big(\Delta u -f\big)v \ d\Omega =0,
\end{equation}
where $v$ is an appropriate test function on $\Omega_\sigma^y$.
Employing the Divergence Theorem, we get an equation
\bql{localweakformeq}
\lambda_{y,\sigma,v}(u):=\int_{\Gamma_\sigma^y\setminus \Gamma_N}(\nabla u\cdot n) v \ d\Gamma - \int_{\Omega_\sigma^y}
\nabla u\cdot \nabla v \ d\Omega =
\int_{\Omega_\sigma^y} fv \  d\Omega-\int_{\Gamma_\sigma^y\cap \Gamma_N}u_N v \ d\Gamma
\end{equation}
of the form (\ref{eqlkubk}).
For nodes whose subdomain boundary does not intersect $\Gamma_N$,
the second term on the right--hand side vanishes.

Note that neither Lagrange multipliers nor penalty parameters are
introduced into the local weak form,
because the Dirichlet boundary
conditions 
are imposed  directly using the second line of
(\ref{eqfullcoll}) for suitable collocation points, usually
taking  a subset of the trial nodes.

Some variations of MLPG differ in their choice of functionals
(\ref{localweakformeq}). If the test function $v$ is chosen to vanish
on $\Gamma_\sigma^y\setminus \Gamma_N$, the first integral
in (\ref{localweakformeq}) is zero, and we have MLPG1.
If the local test function $v$ is the constant 1, the second integral vanishes,
and we have MLPG5.
\section{Implementation}
In this section, we describe the implementation
of GMLS approximations to solve the
Poisson problem (\ref{poissoneq})
using the weak form equations (\ref{localweakformeq}).

At first we fix $m$, the maximal degree of polynomial basis functions we use.
These form the space $\calp:=\mathbb{P}_m^d$ of $d$--variate real--valued polynomials of
degree at most $m$. The dimension of this space is $Q={m+d  \choose d}$.
If the problem
has enough smoothness, $m$ will
determine the convergence rate.

Then we choose a set $X=\{x_1,x_2,...,x_N\}\subset\Omega$ of scattered trial points which
fills 
the domain reasonably well, without letting two points come
too 
close to each other.
To make this more precise, we need the quantities
{\em fill distance} and {\em separation distance} which
are important to measure the quality of centers and
derive rates of convergence. For a set of points $X=\{x_1,x_2,...,x_N\}$ in a bounded domain
$\Omega\subseteq \R^d$, the {\em fill distance} is defined to be
\begin{equation*}
h_{X,\Omega}=\sup_{x\in\Omega}\min_{1\leq j\leq N}\|x-x_j\|_2,
\end{equation*}
and the {\em separation distance} is defined by
\begin{equation*}
q_{X}=\frac{1}{2}\min_{i\neq j}\|x_i-x_j\|_2.
\end{equation*}
A set $X$ of data sites is said to be
{\em quasi-uniform} with respect to a constant $c_{\mathrm{qu}}>0$ if
\bql{quasi-uniform}
q_X\leq h_{X,\Omega}\leq  c_{\mathrm{qu}} q_X.
\end{equation}
In this sense, we require the set $X$ of trial nodes to be quasi--uniform.

We now have to define the functionals
$\lambda_1,\ldots,\lambda_M$ discretizing our PDE problem. This requires a
selection between MLPG1, MLPG2, and MLPG5, and the decision
to use oversampling or not, i.e. $M>N$ or $M=N$.
Oversampling will often increase stability at increased cost,
but we found that in our examples it was not necessary.
Since we have to execute the GMLS method for each functional $\lambda_k$,
approximating it in terms of function values at the trial nodes
in $B(y_k,\delta)\cap X$,
we have to make sure that the GMLS does not break down.
This means that the matrix $\calb$ of (\ref{aWl}) must have  rank $Q$, if
formed for the nodes in $B(y_k,\delta)\cap X$. In general:
\begin{Definition}
A set $Z$ of pairwise distinct points in $\R^d$  is
called $\mathbb P_m^d$-{\em unisolvent} if the zero polynomial is the
only polynomial from $\mathbb P_m^d$ which vanishes on $Z$.
\end{Definition}
To give a sufficient condition for unisolvency, we need
\begin{Definition}
A set $\Omega\subset \mathbb R^d$ is said to
satisfy an {\em interior cone condition} if there
exists an angle $\theta\in(0,\pi/2)$ and a radius $r>0$
such that for every $x\in\Omega$ a unit vector $\xi(x)$ exists such that the cone
$$
C(x,\xi,\theta,r):=\big\{ x+ty: y\in\mathbb R^d, \|y\|_2=1, y^T\xi\geq \cos\theta, t\in[0,r]\big\}
$$
is contained in $\Omega$.
\end{Definition}
\begin{Theorem}\RSlabel{Theunisol} {\rm (\cite{schaback:2011-1}, see also \cite{wendland:2005-1})}\\
For any compact domain $\Omega$ in $\R^d$ with an interior cone condition,
and any $m\geq 0$ there are positive constants $h_0$ and $c_0$  such that
for all trial node sets $X$  with fill distance $h_{X,\Omega}\leq h_0$, all test points
$y\in\Omega$, and all radii $\delta\geq c_0\,h_{X,\Omega}$,
the set $B(y,\delta)\cap X$ is $\mathbb{P}_m^d$--unisolvent.
\end{Theorem}
This means that the placement of test nodes and the choice of weight function
supports can be linked to the fill distance of the set of trial
nodes. Oversampling never causes problems, if the weight function support radius
is kept proportional to the fill distance of the trial nodes.

Some test nodes should be scattered over the Dirichlet
boundary $\Gamma_D$ to impose the Dirichlet
boundary conditions. Like in the collocation case, we 
denote the subset of these points by
$Y_D\subset Y\cap \Gamma_D$. For MLPG2, we similarly define $Y_N$, and then the
setup of the functionals simply follows (\ref{eqfullcoll}), with or without oversampling.
In principle, the sets $Y_\Omega,\,Y_N,\,Y_D$ need not be disjoint.

For weak problems in MLPG1 or MLPG5 form, we just implement the functionals
 $\lambda_{y_k,\sigma_k,v_k}$ of
(\ref{localweakformeq}) as described in Section \ref{Sectweakform}.
Altogether, we follow Section \ref{SecMM} by implementing (\ref{eqlkubk})
via (\ref{eqlkajkuj}), and ending with the system (\ref{eqAub}).
\biglf
The order of convergence of the approximated
functional to its exact value is important in this case.
Applying the same strategy presented in \cite{mirzaei-et-al:2012-1}
 for $\lambda_{y,\alpha}(u):=D^\alpha u(y)$, we can prove
\begin{Theorem}\RSlabel{gmls-error-bound}
Let
$$
\lambda(u)=\lambda_{y,\sigma,w}(u):=\int_{K} w(x) \,Lu(x)dx,
\quad K=\Omega_\sigma^y \mbox{ or } \partial \Omega_\sigma^y,\, y\in\Omega_\sigma^y.
$$
In the
situation of Theorem \ref{Theunisol}, define $\Omega^*$ to be the
closure of $\bigcup_{x\in\Omega}B(x,\delta)$.
Define
$$
\widehat{\lambda(u)}:= \sum_{j=1}^N  a^*_j(\lambda) u(x_j ),
$$
where $a^*_j(\lambda)$ are functions derived from
the GMLS approximation in (\ref{aWl}). Then there exists a
constant $c > 0$ such that for all $u\in C^{m+1}(\Omega^*)$ and all quasi-uniform
$X\subset\Omega$ with $h_{X,\Omega}\leq h_0$ we have
\bql{gmls-error}
\left\|\lambda( u)-\widehat{\lambda(u)}\right\|_{L_\infty(\Omega)}
\leq ch_{X,\Omega}^{m+1-n}|u|_{C^{m+1}(\Omega^*)},
\end{equation}
providing $\int_{K}|w(x)|dx <\infty$ and if $\lambda(u)\neq 0$,
$\int_{K} w(x) \,Lx^\alpha dx\neq 0$ ($\lambda(x^\alpha)\neq 0$)
for some $\alpha$ with $|\alpha|=m$. Here $n$ is the maximal order
of derivatives involved in linear operator $L$ and
$|u|_{C^{m+1}(\Omega^*)}:=\max_{|\alpha|=m+1}\|D^\alpha u\|_{L_\infty(\Omega)}$.
\end{Theorem}


However, we cannot guarantee that the system (\ref{eqAub}) has full rank,
since we only made sure that the rows of the system can be calculated
via the GMLS if Theorem
\ref{Theunisol} applies. Oversampling will usually help if the system
causes problems.

After the solution vector $\u$ of (\ref{eqAub}), consisting of values $u(x_j)$
of values at the trial nodes is determined by solving the system,
other values of the solution function $u(x)$
(and also its derivatives) can be calculated in every point $x\in\Omega$ again using the MLS approximation.

Since we have direct approximations for boundary conditions and
local weak forms, this method is called {\em direct meshless local
Petrov-Galerkin (DMLPG) method}. It comes in the DMLPG1, DMLPG2, and DMLPG5 variations.


In contrast to MLPG2,
if the GMLS derivatives (``diffuse'' derivatives)
\cite{mirzaei-et-al:2012-1} are used instead of the standard
derivatives of MLS shape functions,
we have DMLPG2. As investigated in \cite{mirzaei-et-al:2012-1},
the accuracies for calculating the matrix $\cala$ of (\ref{eqAub})
are the same, but the computational cost of DMLPG2 is less.
When looking into the literature, we found that DMLPG2
coincides with the {\em Diffuse Approximate Method} (DAM) \cite{prax-et-al:1996-1}.
But since we avoid using the word {\em diffuse}
because there is nothing ``diffuse'' about these derivatives \cite{mirzaei-et-al:2012-1},
we will call the method DMLPG2 or {\em direct MLS collocation (DMLSC) method}.

As we saw in Section \ref{SecMLS}, in DMLPG1/5 methods
the integrations are done only over polynomials, if the latter are used in the
GMLS. For every
functional $\lambda_k$, $1\leq k\leq M\geq N$, the GMLS routine is
called only once. There are no calls to produce values of shape functions.
The standard MLPG/MLS  technique at degree $m$
implements numerical integration by calling shape function evaluations,
and thus the MLS routine is called approximatively $M\,K$ times where $K$ is the
average number of integration
points. Moreover, in standard MLPG methods the derivatives of MLS shape
function must also be provided, while DMLPG has no shape functions at all.
Consequently, DMLPG is considerably faster than MLPG. In addition,
due to the error analysis presented in
Theorem \ref{gmls-error-bound}
for the new GMLS method, the final accuracies of
both MLPG and DMLPG methods are the same.
We will see experimentally that DMLPG is even more accurate than MLPG.

As highlighted in \cite{babuska-et-al:2009-1}, numerical integration in
FEM is simple because the integrands of the elements of the
stiffness matrix are polynomials. In contrast to this, the shape functions
used in standard meshless methods are much more costly to evaluate,
making numerical integration a much bigger challenge than for the FEM.
In MLPG methods, numerical integrations are simpler than for
various other meshless
methods, since the local weak form breaks everything down to
local well--shaped subdomains.
However, since the integrands are based on MLS shape
functions and their derivatives, a Gauss quadrature with many points is
required to get accurate results, especially when the density of nodes
increases.
Overcoming this drawback is a major advantage of DMLPG methods, because the
integrations are done over polynomials, like in FEM.

It is interesting to
note that if local sub-domains are chosen in DMLPG5  as $S(x,\sigma)$
(square or cube), a $(d-1)$--times $\left\lceil\frac{m}{2}\right\rceil$--points
Gauss quadrature gives the exact solution for local boundary
integrals around the nodes in the interior of $\Omega$.
In DMLPG1, if again $S(x,\sigma)$ is chosen as a local sub-domain
and if a polynomial test function is employed, a $d$-times
$\left\lceil\frac{(m-1)(n-1)+1}{2}\right\rceil$--points Gauss quadrature is
enough to get exact interior local domain integrals.
Here, $n$ is the degree of the polynomial test function.
As a polynomial test function on the square or cube for DMLPG1 with $n=2$, we can use
\begin{equation*}\RSlabel{testfunction}
v(x;x_k)=\left\{
\begin{split}
\prod_{i=1}^d\left(1-\frac{4}{\sigma^2}
(\xi_i-\xi_{ki})^2\right) ,&\quad x\in S(x_k,\sigma), \\
0,\quad\quad\quad\quad\quad\quad\quad\quad & \quad  \mbox{otherwise}
\end{split}\right.
\end{equation*}
where $x=(\xi_1, ... ,\xi_d)$ and
$x_k=(\xi_{k1},...,\xi_{kd})$.
In DMLPG1 with balls as sub-domains, weight functions of the form
function 
\bql{weight}
w_\delta(x,y)=\phi\left( \frac{\|x-y\|_2}{\delta}\right),
\end{equation}
with $\delta = \sigma$
can be used as test functions.
Both of these test functions vanish on $\Gamma_\sigma^{x_k}\backslash\Gamma_N $, as required in DMLPG1.

Note that, if the second weak forms (Green forms) are
taken over local sub-domains and a modified fundamental solution is
used as test function,
the process gives the DMLPG4 rather than MLPG4 or
the meshless LBIE method presented in \cite{zhu-et-al:1998-1}.
In DMLPG4, it is better to use balls as local sub-domains,
because in this case the modified fundamental solution, used as a test function,
can be derived easily. But the test function is not a polynomial.

{The trial and test functions in both MLPG3 and 6 come from the same space and thus they
are Galerkin type techniques. If we formulate the analogues DMLPG3/6, the integrands include shape functions again. Therefore, they annihilate
the advantages of DMLPG methods with respect to numerical integration, and we ignore them in favour of keeping all benefits of DMLPG
methods.}

Instead, we add some remarks about selecting $m$,
the degree of polynomial basis functions
in the GMLS. For $m=1$, the variants DMLPG 1, 4, and 5 will
necessarily fail.
The background is that the GMLS performs an optimal recovery
of a functional $\lambda$ in terms of nodal values, and the recovery is exact
on a subspace $\mathcal{P}$, using minimal coefficients at the nodal values.
Thus, in all cases where the functional is zero on $\mathcal{P}$
by some reason or other,
the recovery formula will be zero and will generate a zero row in
the stiffness matrix.
This happens for all variations based on functionals (\ref{localweakformeq})
and functionals extracted from the second weak form
on interior points,
since all those functionals are reformulations of
$$
\lambda_{y,\sigma,v}(u)=\int_{B(y,\sigma)}v\, \Delta u\,d\omega
$$
and thus vanish on harmonic functions $u$, in particular on linear functions.
Thus, for solving inhomogeneous problems,
users should pick spaces $\calp$ of non--harmonic functions,
if they employ GMLS with exactness on $\calp$. This rules out polynomials with
degree $m\leq 1$.

Another closely related point arises from symmetry of subdomains.
Since polynomials in a ball $B(x,\sigma)$ or a cube $S(x,\sigma)$
have symmetry properties, the entries of stiffness matrices
in rows corresponding to internal points will often be the same
for $m=2k$ and $m=2k+1$. Thus convergence rates often do not increase
when going from $m=2$ to $m=3$, for instance. But this observation affects
MLPG and DMLPG in the same way.
\section{Stability and convergence}
For the classical MLS and the generalized MLS from
\cite{mirzaei-et-al:2012-1}
and Theorem \ref{gmls-error-bound}
it is known that the
recovery
$\widehat {\lambda (u^*)}$
 of values of functionals $\lambda$ on a true solution $u^*$ has an error
of order $\mathcal{O}(h^{m+1-k})$, if $h$ is the fill distance of the
trial nodes, $m$ is the degree of polynomials used locally,
if the exact solution $u^*$ is at least $C^{m+1}$,
$k$ is the maximal order of derivatives of $u^*$ involved in the functional,
and if numerical integration has an even smaller error.
In particular, the classical MLS and the new GMLS produce roughly the same
stiffness matrices, but the GMLS has a considerably smaller computational complexity.



However, the error committed in the approximation of the test functionals
in terms of function values at nodes does not always carry over to
the convergence rate of the full algorithm, since there is no stability
analysis, so far. Even if perfect stability would hold, the best one can expect
is to get the convergence rate implied by the local trial approximation,
i.e. by local polynomials of degree $m$. This would again mean a rate of
$\mathcal{O}(h^{m+1-k})$, but only if the solution is indeed locally
approximated by polynomials of that degree. In fact, the next section will show
that this rate can often be observed. But our symmetry arguments
at the end of the previous section show that sometimes the degree $m=2k+1$
cannot improve the behavior for $m=2k$, because the odd--degree polynomials
simply do not show up in most of the calculations for the stiffness matrix.
\section{Numerical results}
In this section some numerical results are presented to
demonstrate the efficiency of DMLPG methods and its advantages over MLPG methods.
We consider the Poisson equation (\ref{poissoneq})
on the square $[0,1]^2\subset \R^2$ with Dirichlet
boundary conditions.
Since we want to study convergence rates without being limited by
smnoothness of the solution,
we take  {\em Franke's function} \cite{franke:1982-1}
\begin{equation*}
\begin{split}
u(\bar x,\bar y)=&\frac{3}{4}e^{-1/4((9\bar x-2)^2
+(9\bar y-2)^2)}+\frac{3}{4}e^{-(1/49)(9\bar x+1)^2-(1/10)(9\bar y+1)^2)}\\
&~+\frac{1}{2}e^{-1/4((9\bar x-7)^2+(9\bar y-3)^2)}-\frac{1}{5}e^{-(9\bar x-4)^2-(9\bar y-7)^2},
\end{split}
\end{equation*}
where $(\bar x,\bar y)$ denotes the two components of $x\in\R^2$,
as analytical solution and calculate the right hand side and
boundary conditions accordingly. Note that Franke's function is a
standard test function for 2D scattered data fitting.
Regular mesh distributions with mesh-size $h$ are provided in all cases,
though the methods would work with scattered data.
We do not implement oversampling in the results of this paper. In fact, the
trial and test points are chosen to be coincident.
Also, the {\em shifted scaled polynomial}
\begin{equation*}
\left\{\frac{(x-z)^\alpha}{h^{|\alpha|}}\right\}_{0\leq|\alpha|\leq m},
\end{equation*}
where $z$ is a fixed evaluation point such as
a test point or a Gaussian point for integration, is used
instead of the natural polynomial basis $\{x^\alpha\}_{0\leq|\alpha|\leq m}$
for MLS approximation. In \cite{mirzaei-et-al:2012-1}, it is
shown that this choice of basis function
leads to more stable and accurate MLS approximation.
We use the shifted
scaled basis for both MLPG5 and DMLPG5 methods with $m=2,3$ and $4$.
The Gaussian weight function
\begin{equation*}
w_\delta(x,x_j)=\left\{
\begin{split}
\frac{\exp\big(-(\|x-x_j\|_2/c)^2\big)-\exp\big(-(\delta/c)^2\big)}{1-\exp\big(-(\delta/c)^2\big)},&\quad 0\leq \|x-x_j\|_2\leq \delta, \\
0,\quad\quad\quad\quad\quad\quad\quad\quad & \quad  \|x-x_j\|_2>\delta
\end{split}\right.
\end{equation*}
is employed where $c=c_0h$ is a constant controlling the
shape of the weight function and $\delta=\delta_0h$ is the size of the support domains.

Let $m=2$ and set $c_0=0.6$ and $\delta_0=2m$. At first
the local sub-domains are taken to be circles. To get the
best results in MLPG we have to use an accurate quadrature
formula. Here a 20-points regular Gauss-Legendre quadrature
is employed for numerical integrations over local sub-domains.

Numerical results, for different mesh-sizes $h$, are presented
in terms of maximum errors, ratios and CPU times used for MLPG5 and DMLPG5 in Tables 1.
\begin{center}
\begin{tabular}{llllllllllllll}
  \multicolumn{14}{l}{\footnotesize{\textsc{Table 1}. The maximum errors, ratios and CPU times used for MLPG5 and DMLPG5 with $m=2$}}\\
  \hline
  & \fs{}       && \fs{MLPG5}     &&\fs{}                 && \fs{DMLPG5}           &&\fs{}         &&\fs{CPU time used} && \fs{}      \\
                \cline{4-6}                               \cline{8-10}                                                  \cline{12-14}
  & \fs$h$      && \fs{$\|e\|_\infty$}      &&\fs{ratio}   && \fs{$\|e\|_\infty$}  &&\fs{ratio}    && \fs{MLPG5}        && \fs{DMLPG5}   \\
  \hline
  & \fs$0.2$    && \fs$0.44\times 10^{-1}$ && \fs$-$      &&\fs$0.23\times10^{-1}$ &&\fs$-$        &&\fs {$0.4$ sec.}   && \fs{$0.2$ sec.}  \\
  & \fs$0.1$    && \fs$0.15\times 10^{-1}$ && \fs$1.59$   &&\fs$0.72\times10^{-2}$ &&\fs$1.68$     &&\fs $1.2$          && \fs$0.3$  \\
  & \fs$0.05$   && \fs$0.73\times 10^{-2}$ && \fs$0.99$   &&\fs$0.20\times10^{-2}$ &&\fs$1.84$     &&\fs $6.5$         && \fs$1.4$  \\
  & \fs$0.025$  && \fs$0.24\times 10^{-2}$ && \fs$1.61$   &&\fs$0.58\times10^{-3}$ &&\fs$1.80$     &&\fs $68.5$        && \fs$6.5$ \\
  & \fs$0.0125$ && \fs$0.66\times 10^{-3}$ && \fs$1.85$   &&\fs$0.14\times10^{-3}$ &&\fs$1.98$     &&\fs $2016.0$       && \fs$52.1$ \\

 \hline
 \newline
 \newline
\end{tabular}
\end{center}
The mesh-size $h$ is divided by two row by row, therefore the ratios are computed by
$$\log_2\left(\frac{\|e({h})\|_\infty}{\|e({h/2})\|_\infty}\right).$$
Both methods have nearly the same order 2, which cannot be improved for this
trial space, since the expected optimal order is $m+1-k=3-1=2$.
But significant
differences occur in the columns with CPU times.
As we stated before, this is due to restricting local integrations
to polynomial basis functions in DMLPG rather than to integrate over MLS shape
functions in the original MLPG. We could get the same results
with fewer integration points for DMLPG, but to be fair in comparison, we use the same quadrature.

In addition, to give more insight into the errors, the maximum
errors of MLPG5 and DMLPG5 are illustrated in Fig. \ref{fig2}.
Once can see that DMLPG is more accurate, maybe due to avoiding
many computations and hence many roundoff errors.

\begin{figure}[hbt]
\begin{center}
\includegraphics[width=10cm]{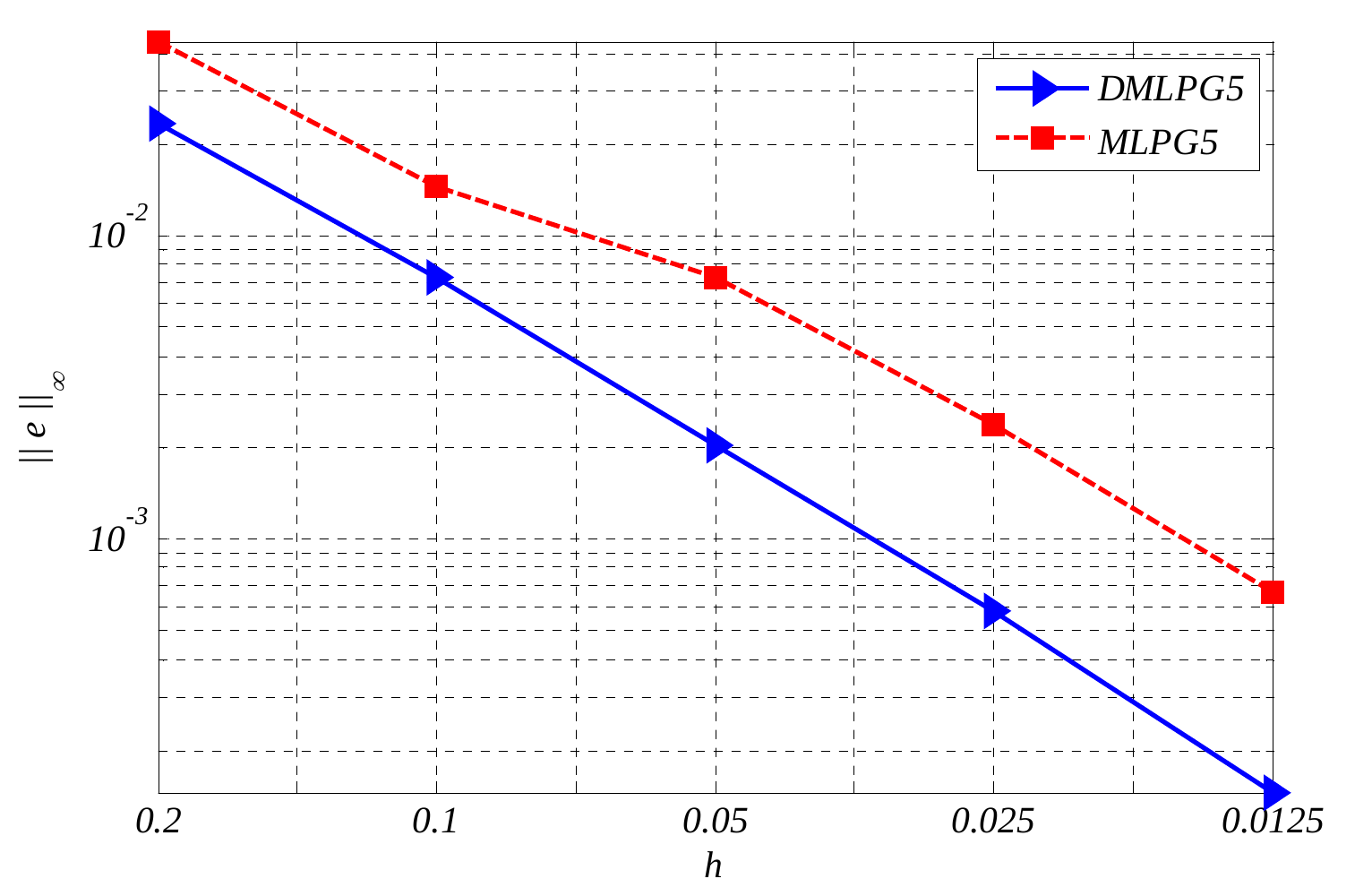}\\
\caption{\small{Comparison of MLPG5 and DMLPG5 in terms of maximum errors for $m=2$.}}\RSlabel{fig2}
\end{center}
\end{figure}
In Table 2 and Fig. \ref{fig3}, we have compared MLPG5 and DMLPG5 in case $m=3$.
The convergence rate stays at 2 for both methods, since the third--degree
polynomials cannot contribute to the weak form. The figure shows approximately the same
error results. But the CPU times used are indeed different.
\begin{center}
\begin{tabular}{llllllllllllll}
  \multicolumn{14}{l}{\footnotesize{\textsc{Table 2}. The maximum errors, ratios and CPU times used for MLPG5 and DMLPG5 with $m=3$}}\\
  \hline
  & \fs{}       && \fs{MLPG5}     &&\fs{}                 && \fs{DMLPG5}           &&\fs{}         &&\fs{CPU time used} && \fs{}      \\
                \cline{4-6}                               \cline{8-10}                                                  \cline{12-14}
  & \fs$h$      && \fs{$\|e\|_\infty$}      &&\fs{ratio}   && \fs{$\|e\|_\infty$}  &&\fs{ratio}    && \fs{MLPG5}        && \fs{DMLPG5}   \\
  \hline
  & \fs$0.2$    && \fs$0.28\times 10^{-1}$ && \fs$-$      &&\fs$0.23\times10^{-1}$ &&\fs$-$        &&\fs {$0.8$ sec.}   && \fs{$0.2$ sec.}  \\
  & \fs$0.1$    && \fs$0.13\times 10^{-1}$ && \fs$1.08$   &&\fs$0.74\times10^{-2}$ &&\fs$1.62$     &&\fs $1.8$         && \fs$0.4$  \\
  & \fs$0.05$   && \fs$0.33\times 10^{-2}$ && \fs$1.98$   &&\fs$0.20\times10^{-2}$ &&\fs$1.89$     &&\fs $9.7$        && \fs$1.6$  \\
  & \fs$0.025$  && \fs$0.78\times 10^{-3}$ && \fs$2.09$   &&\fs$0.58\times10^{-3}$ &&\fs$1.80$     &&\fs $87.7$        && \fs$7.6$ \\
  & \fs$0.0125$ && \fs$0.19\times 10^{-3}$ && \fs$2.06$   &&\fs$0.15\times10^{-3}$ &&\fs$1.98$     &&\fs $2293.3$       && \fs$56.1$ \\

 \hline
 \newline
 \newline
\end{tabular}
\end{center}
\begin{figure}
\begin{center}
\includegraphics[width=10cm]{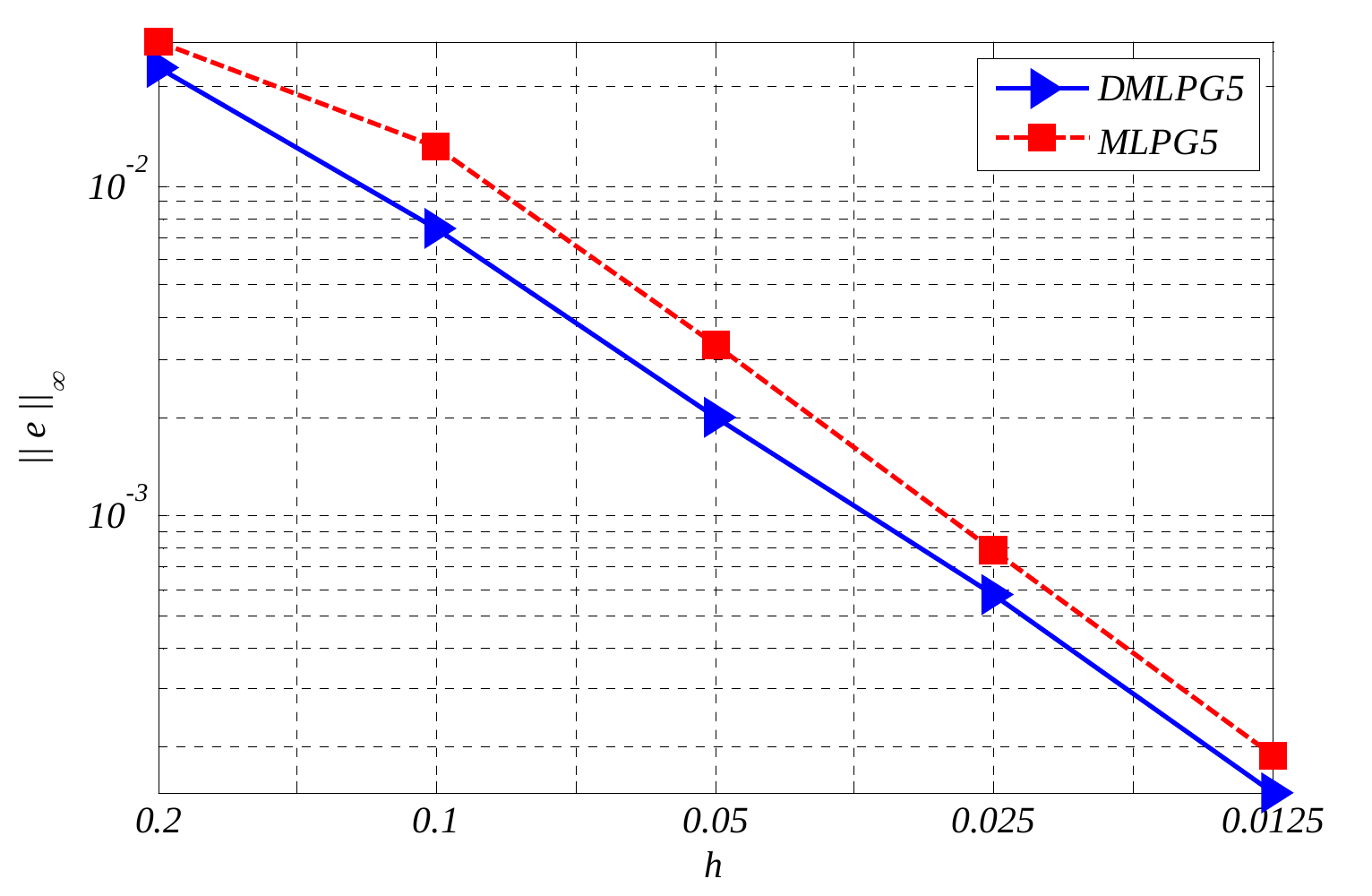}\\
\caption{\small{Comparison of MLPG5 and DMLPG5 in terms of maximum errors for $m=3$.}}\RSlabel{fig3}
\end{center}
\end{figure}


In the results presented up to here, we used balls (circles) as local sub-domains.
Now we turn to use squares for both MLPG5 and DMLPG5.
Also, we run the programs with $m=4$ to see the differences.
The parameters $c_0=0.8$ and $\delta_0=2m$ are selected.
In DMLPG5, a 2-points Gaussian quadrature is enough to get exact numerical
integrations.
But for MLPG5 and the right hand sides
we use a 10-points Gaussian quadrature for every sides of squares.
The results are depicted in Table 3 and Fig. \ref{fig4}.
DMLPG is more accurate and approximately gives the full order $4$
in this case. Beside, as we expected, the computational
cost of DMLPG is remarkably less than MLPG.
\begin{center}
\begin{tabular}{llllllllllllll}
  \multicolumn{14}{l}{\footnotesize{\textsc{Table 3}. The maximum errors, ratios and CPU times used for MLPG5 and DMLPG5 with $m=4$}}\\
  \hline
  & \fs{}       && \fs{MLPG5}              &&\fs{}        && \fs{DMLPG5}           &&\fs{}         &&\fs{CPU time used} && \fs{}      \\
                \cline{4-6}                               \cline{8-10}                                                  \cline{12-14}
  & \fs$h$      && \fs{$\|e\|_\infty$}     &&\fs{ratio}   && \fs{$\|e\|_\infty$}   &&\fs{ratio}    && \fs{MLPG5}        && \fs{DMLPG5}   \\
  \hline
  & \fs$0.2$    && \fs$0.10\times 10^{0}$  && \fs$-$      &&\fs$0.12\times10^{0}$  &&\fs$-$        &&\fs {$0.5$ sec.}   && \fs{$0.2$ sec.}  \\
  & \fs$0.1$    && \fs$0.25\times 10^{-1}$ && \fs$2.04$   &&\fs$0.17\times10^{-1}$ &&\fs$2.87$     &&\fs $2.7$         && \fs$0.2$  \\
  & \fs$0.05$   && \fs$0.78\times 10^{-2}$ && \fs$1.66$   &&\fs$0.12\times10^{-2}$ &&\fs$3.75$     &&\fs $19.2$        && \fs$0.9$  \\
  & \fs$0.025$  && \fs$0.79\times 10^{-3}$ && \fs$3.30$   &&\fs$0.75\times10^{-4}$ &&\fs$4.04$     &&\fs $142.2$       && \fs$4.7$ \\
  & \fs$0.0125$ && \fs$0.55\times 10^{-4}$ && \fs$3.86$   &&\fs$0.43\times10^{-5}$ &&\fs$4.12$     &&\fs $2604.9$       && \fs$43.9$ \\

 \hline
 \newline
\end{tabular}
\end{center}
\begin{figure}
\begin{center}
\includegraphics[width=10cm]{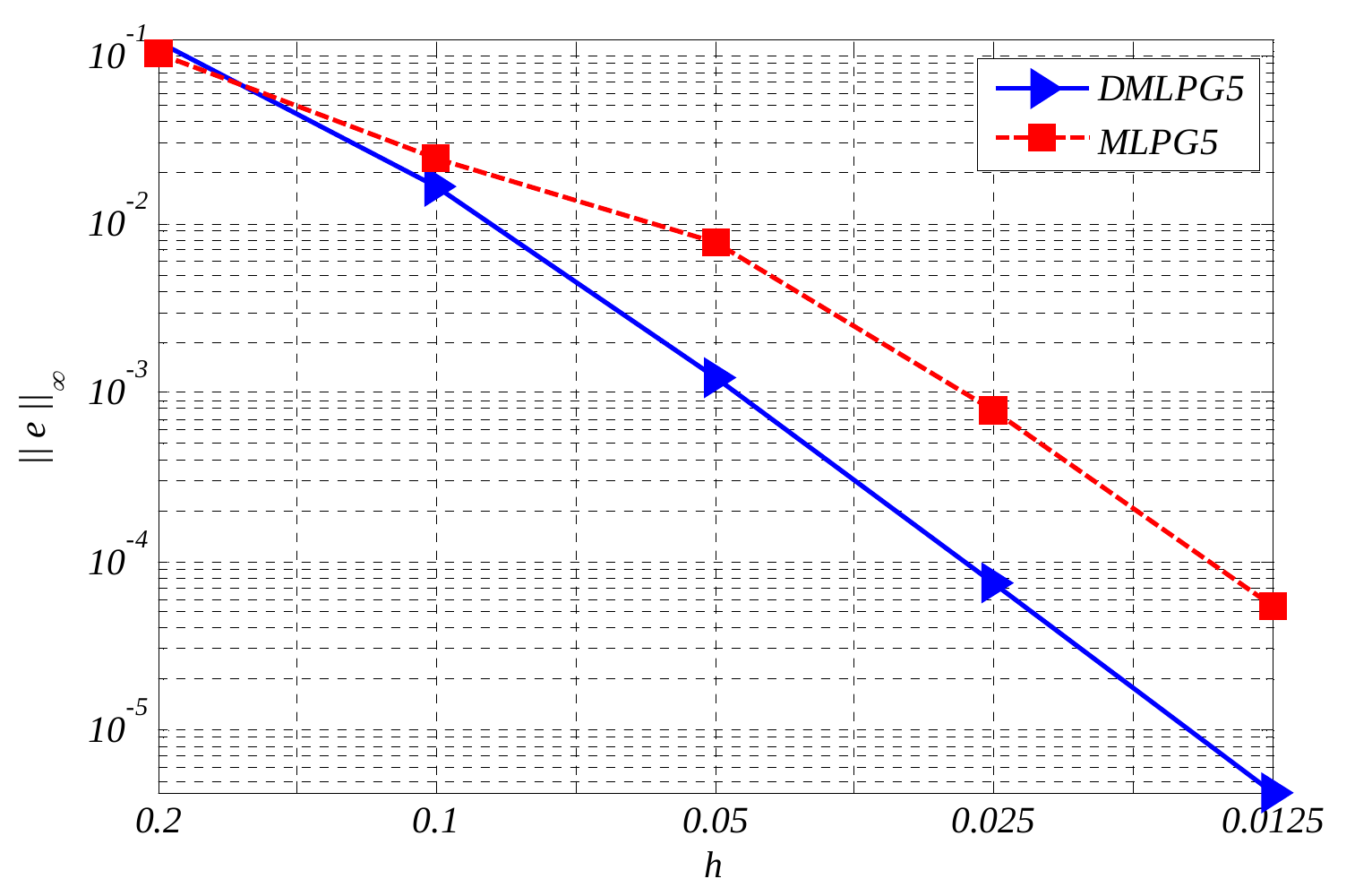}\\
\caption{\small{Comparison of MLPG5 and DMLPG5 in terms of maximum errors for $m=4$.}}\RSlabel{fig4}
\end{center}
\end{figure}


Results for MLPG1 and DMLPG1 turn out to behave similarly.
As we know, MLPG1 is more expensive than MLPG5 \cite{atluri:2005-1,atluri-shen:2002-book},
but there is no significant difference between computational
costs of DMLPG5 and DMLPG1. Therefore the differences between
CPU times used for MLPG1 and DMLPG1 are absolutely larger.

All routines were written using \textsc{Matlab}$^\copyright$
 and run on a Pentium 4 PC with 2.50 GB of Memory and a twin--core 2.00 GHz CPU.

\section{Conclusion}
This article describes a new MLPG method,
called DMLPG method, based on generalized moving least squares (GMLS)
approximation for solving boundary value problems. The new method is
considerably faster than the classical MLPG variants, because
\begin{itemize}
\item direct approximations of data functionals are used
for Dirichlet boundary conditions and local weak forms,
\item local integrations are done over
polynomials rather than over complicated MLS shape
functions,
\item numerical integrations can sometimes be performed exactly.
\end{itemize}
The convergence rate of both methods should be the same,
but thanks to avoiding many computations and roundoff errors,
and of course by treating the numerical integrations in a more elegant and
possibly exact
way, the results of DMLPG turn often out to be more accurate than the results of MLPG.

On the downside, DMLPG does not work for $m=1$ since it locally uses (harmonic) linear functions
instead of complicated shape functions. But most  MLPG users
choose higher degrees anyway, in order to obtain better
convergence rates.

Altogether, we believe that the DMLPG methods have great potential to
replace the original MLPG methods in many situations.
\section*{Acknowledgments}
The first author was financially supported by the Center of Excellence for Mathematics, University of Isfahan. 
\bibliographystyle{elsart-num-sort}


\end{document}